\documentclass[aip,amsmath,amssymb,reprint]{revtex4-1}

\usepackage{dcolumn}
\usepackage{amsmath}
\usepackage{txfonts}
\usepackage{cancel}
\usepackage[T1]{fontenc}
\usepackage{mathptmx}

\begin{document}

\title{Sample variance of rounded variables}

\author{J.~An}
\email{0x1b6f@gmail.com}
\affiliation{$^1$ Korea Astronomy and Space Science Institute}

\date{\today}

\begin{abstract}
If the rounding errors are assumed to be distributed independently
from the intrinsic distribution of the random variable,
the sample variance $s^2$ of the rounded variable is given by the sum
of the true variance $\sigma^2$ and the variance of the rounding errors
(which is equal to $w^2/12$ where $w$ is the size of the rounding window).
Here the exact expressions for the sample variance of the rounded variables
are examined and it is also discussed when the simple approximation
$s^2=\sigma^2+w^2/12$ can be considered valid. In particular,
if the underlying distribution $f$ belongs to a family of symmetric
normalizable distributions such that $f(x)=\sigma^{-1}F(u)$ where
$u=(x-\mu)/\sigma$, and $\mu$ and $\sigma^2$ are the mean and variance
of the distribution, then the rounded sample variance scales
like $s^2-(\sigma^2+w^2/12)\sim\sigma\Phi'(\sigma)$ as $\sigma\to\infty$
where $\Phi(\tau)=\int_{-\infty}^\infty{\rm d}u\,e^{iu\tau}F(u)$
is the characteristic function of $F(u)$. It follows that, roughly speaking,
the approximation is valid for a slowly-varying symmetric underlying
distribution with its variance sufficiently larger than the size of
the rounding unit.
\end{abstract}

\maketitle

\section{Introduction}

Most real world data are only recorded in the rounded figure with a fixed
number of significant digits. Strictly this rounding introduces additional
systematic uncertainties which must be properly accounted for,
in order to infer the property of the intrinsic distribution of
the measured quantities. Naively, assuming that there is neither
intrinsic uncertainty nor systematic bias, the differences between
the true value and the rounded reported value are expected to be
distributed evenly over the window of the size of the reporting unit.

In particular, if the variable value $x$ is rounded to an integer multiple
value of the measurement unit as in $nw$ (where $w$ is the measurement
unit and $n$ is an integer), then
$(n+\delta-1/2)w\le x<(n+\delta+1/2)w$ or
$(n+\delta-1/2)w<x\le (n+\delta+1/2)w$. Here the constant
$\delta\in[-1/2,1/2]$ specifies the rounding method
(e.g., $\delta=1/2$ for rounding down to the floor,
$\delta=-1/2$ for rounding up to the ceiling, and
$\delta=0$ for rounding to the nearest integer etc.)
whereas the equal signs at the boundary follow the prescribed convention.
Then the rounding error (i.e.\ $\rho=nw-x$) is distributed
in the rectangular distribution:
\begin{equation}
P(\rho)=\begin{cases}w^{-1}&\text{for $-(1/2+\delta)w<\rho<(1/2-\delta)w$}
\\0&\text{elsewhere}\end{cases},
\end{equation}
where the distribution at the boundary is determined by the chosen convention
-- however provided that $x$ is a real variable in a continuous distribution,
the boundaries consitute a null measure set and so the specific choice does not
affect the following discussion. For a random variable $x$,
the mean of the rounded values is (with $\bar x$ being the true mean of $x$)
\begin{equation}
\overline{nw}=\overline{x+\rho}=\bar x
+\int_{-(1/2+\delta)w}^{(1/2-\delta)w}\frac\rho w\,{\rm d}\rho
=\bar x-\delta u,
\end{equation}
while the variance is
\begin{equation}\begin{split}
s^2&=\overline{(nw-\overline{nw})^2}
=\overline{(x+\rho)^2}-(\bar x+\bar\rho)^2
\\&=\sigma^2+\overline{\rho^2}-\bar\rho^2+2(\overline{x\rho}-\bar x\bar\rho),
\end{split}\end{equation}
where $\sigma^2=\overline{x^2}-\bar x^2$ is the variance of $x$,
with the variance of the rounding errors given by
\begin{equation}
\overline{\rho^2}-\bar\rho^2
=\int_{-(1/2+\delta)w}^{(1/2-\delta)w}\frac{\rho^2}w\,{\rm d}\rho-(\delta w)^2
=\frac{w^2}{12}.
\end{equation}
In other words, provided that the distribution of $x$ does not affect
the rounding (as in $\overline{x\rho}=\bar x\bar\rho$),
the standard deviation of the rounded values is simply
a quadrature sum of the true underlying standard deviation
and that of the rounding errors, and the true standard deviation
may be estimated from the variance of the rounded values via
\begin{equation}
\sigma=\left(s^2-\frac{w^2}{12}\right)^{1/2}.
\end{equation}
However, this result is only valid ``on average'' sense. That is to say,
the underlying distribution of the variable can technically affect
the rounding but for an arbitrary unspecified distribution,
the expected value of ``$\overline{x\rho}-\bar x\bar\rho$'' should be zero
and the reported error tends to the quadrature sum of the random
error and the rounding error ($\sigma_\rho/w=1/\!\sqrt{12}\approx0.2887$).

\section{Theory}

Suppose that $f(x)$ is a probability distribution of a real random
variable $x$ with
\begin{equation}\label{eq:norm}\begin{split}
&\int_{-\infty}^\infty\!{\rm d}x\,f(x)=1
,\quad
\int_{-\infty}^\infty\!{\rm d}x\,x\,f(x)=\mu,\\&
\int_{-\infty}^\infty\!{\rm d}x\,(x-\mu)^2f(x)=
\int_{-\infty}^\infty\!{\rm d}x\,x^2f(x)-\mu^2=\sigma^2.
\end{split}\end{equation}
Next consider the rounding off the measured value of the variable
such that, with a fixed constant $\delta\in[-1/2,1/2]$ and the measurement
unit $w$, the value $x$ is read off by an integer multiple of the unit,
i.e.\ $nw$, where $n=\left\lfloor(x/w+1/2-\delta)\right\rfloor$ or 
$n=\left\lceil(x/w-1/2-\delta)\right\rceil$ with $\lfloor x\rfloor$
and $\lceil x\rceil$ being the integer floor and ceiling of $x$.
Then the (discrete) distribution of the reported integer $n$ for
the rounded value is found to be
\begin{equation}\label{eq:fndef}
F_n=\int_{(n+\delta-1/2)w}^{(n+\delta+1/2)w}\!{\rm d}x\,f(x).
\end{equation}
This distribution is properly normalized: that is,
\begin{equation}
\sum_{n=-\infty}^\infty F_n=\int_{-\infty}^\infty\!{\rm d}x\,f(x)=1,
\end{equation}
and so we can find the mean and the variance of the rounded variables
by calculating
\begin{equation}\begin{split}
m&=\bar nw=w\sum_{n=-\infty}^\infty nF_n,\quad
\\s^2&=\overline{(nw-m)^2}=\overline{n^2}w^2-m^2
\\&=w^2\left[\sum_{n=-\infty}^\infty n^2F_n
-\left(\sum_{n=-\infty}^\infty nF_n\right)^2\right].
\end{split}\end{equation}
For some distributions $f(x)$, the associated discrete distribution $F_n$
as well as its mean $m$ and the standard deviation $s^2$ of
the rounded variable can be computed analytically. However the calculations
become quite tedious even for many simple distributions and the computations
can only be done numerically for most distributions including the important
example such as the normal distribution. 
Instead here we try to analyze the problem more generally.
Henceforth we also assume $w=1$ but the requisite adjustments for
any other value of $w$ are trivial.

\subsection{characteristic function}

First let us introduce the characteristic function $\varphi(t)$
of the distribution $f(x)$: namely,
\begin{equation}
\varphi(t)=\int_{-\infty}^\infty\!{\rm d}x\,e^{itx}f(x).
\end{equation}
The derivatives of the characteristic function then result in
\begin{equation}\begin{split}
\varphi^{(n)}(t)&=i^n\int_{-\infty}^\infty\!{\rm d}x\,x^ne^{itx}f(x);
\\\varphi^{(n)}(0)&=i^n\int_{-\infty}^\infty\!{\rm d}x\,x^nf(x),
\end{split}\end{equation}
and so $\varphi(0)=1$, $\varphi'(0)=i\mu$ and
$\varphi''(0)=-(\sigma^2+\mu^2)$.
We can also define the shifted characteristic function:
\begin{equation}\label{eq:scfdef}\begin{split}
\tilde\varphi(t)&=e^{-it\mu}\varphi(t)
=\int_{-\infty}^\infty\!{\rm d}x\,e^{it(x-\mu)}f(x)
\\&=\int_{-\infty}^\infty\!{\rm d}\epsilon\,e^{it\epsilon}f(\mu+\epsilon),
\\\tilde\varphi'(t)&=e^{-it\mu}[\varphi'(t)-i\mu\varphi(t)];
\\\tilde\varphi''(t)&=e^{-it\mu}
[\varphi''(t)-2i\mu\varphi'(t)-\mu^2\varphi(t)].
\end{split}\end{equation}
Then $\tilde\varphi(0)=\varphi(0)=1$,
$\tilde\varphi'(0)=\varphi'(0)-i\mu\varphi(0)=0$, and
$\tilde\varphi''(0)=\varphi''(0)-2i\mu\varphi'(0)-\mu^2\varphi(0)
=-\sigma^2-\mu^2+2\mu^2-\mu^2=-\sigma^2$.
In other words, the Maclaurin series coefficients of $\tilde\varphi(t)$
result in the sequence of the central moments whereas those of
$\varphi(t)$ result in the moments about the origin. Furthermore,
if the distribution is symmetric about its mean $\mu$ as in
$f(\mu+\epsilon)=f(\mu-\epsilon)$ for any $\epsilon\in\mathbb R$, then
\begin{equation}\begin{split}
\tilde\varphi(t)
&=\int_{-\infty}^\infty\!{\rm d}\epsilon\,e^{it\epsilon}f(\mu-\epsilon)
\\&=\int_{-\infty}^\infty\!{\rm d}\varepsilon\,
e^{i(-t)\varepsilon}f(\mu+\varepsilon)=\tilde\varphi(-t);
\end{split}\end{equation}
and also $\tilde\varphi^{(n)}(t)=(-1)^n\tilde\varphi^{(n)}(-t)$.
That is to say, the shifted characteristic function of a symmetric
distribution is an even function. The converse also holds in that,
if the characteristic function is in the form of
$\varphi(t)=e^{it\mu}\tilde\varphi(t)$ with an even function such that
$\tilde\varphi(t)=\tilde\varphi(-t)$, the distribution must be
symmetric about the mean $\mu$.

\subsection{distribution of rounded values}

The characteristic function may also be inverted
to recover the distribution via the inverse Fourier transform: that is,
\begin{equation}\label{eq:invf}\begin{split}
f(x)&=\frac1{2\pi}\int_{-\infty}^\infty\!{\rm d}t\,e^{-itx}\varphi(t)
\\&=\frac1{2\pi}\int_{-\infty}^\infty\!{\rm d}t\,e^{it(\mu-x)}\tilde\varphi(t).
\end{split}\end{equation}
Inserting this into equation (\ref{eq:fndef}), we find
the expression for the discrete distribution $F_n$ of the
rounded variable in terms
of the characteristic function $\varphi(t)$: namely,
\begin{equation}\begin{split}
F_n&=\frac1{2\pi}\int_{n+\delta-1/2}^{n+\delta+1/2}\!{\rm d}x
\int_{-\infty}^\infty\!{\rm d}t\,e^{-itx}\varphi(t)
\\&=\frac1{2\pi}\int_{-\infty}^\infty\!{\rm d}t\,\varphi(t)
\int_{n+\delta-1/2}^{n+\delta+1/2}\!{\rm d}x\,e^{-itx}
\\&=\frac1{2\pi}\int_{-\infty}^\infty\!{\rm d}t\,\varphi(t)
\,\mbox{sinc}\Bigl(\frac t2\Bigr)\,e^{-it(n+\delta)}
\\&=\frac1{2\pi}\int_{-\infty}^\infty\!{\rm d}t\,\tilde\varphi(t)
\,\mbox{sinc}\Bigl(\frac t2\Bigr)\,e^{it(\mu-\delta-n)}.
\end{split}\end{equation}
Here $\mbox{sinc}(x)=x^{-1}\sin x$ for $x\ne0$ and $\mbox{sinc}(0)=1$.
In addition we can also define the characteristic function of $F_n$.
Since $F_n$ is a discrete distribution, its characteristic
function is given by\citep{Tr84}
\begin{equation}\label{eq:fn_char}\begin{split}
\Phi_t&=\sum_{n=-\infty}^\infty\!e^{itn}F_n
\\&=\frac1{2\pi}\int_{-\infty}^\infty\!{\rm d}\tau\,\varphi(\tau)
\,\mbox{sinc}\Bigl(\frac\tau2\Bigr)\,e^{-i\tau\delta}
\sum_{n=-\infty}^\infty\!e^{i(t-\tau)n}
\\&=\sum_{k=-\infty}^\infty\!\varphi(t+2\pi k)
\,\mbox{sinc}\left(\frac t2+\pi k\right)\,e^{-i(t+2\pi k)\delta}
\\&=\sum_{k=-\infty}^\infty\!\tilde\varphi(t+2\pi k)
\,\mbox{sinc}\left(\frac t2+\pi k\right)\,e^{i(t+2\pi k)(\mu-\delta)},
\end{split}\end{equation}
Here we have used the Fourier series representation of the so-called
Dirac comb distribution: namely,
\begin{equation}
\frac1{2\pi}\sum_{n=-\infty}^\infty\!e^{in(t-\tau)}
=\sum_{k=-\infty}^\infty\!\deltaup(t-\tau+2\pi k).
\end{equation}
Then the derivative of $\Phi_t$ is found to be
\begin{multline}
\frac{d\Phi_t}{dt}
=\sum_{k=-\infty}^\infty\Biggl\{
\left[\varphi'(t+2\pi k)-i\delta\varphi(t+2\pi k)\right]
\,\mbox{sinc}\Bigl(\frac t2+\pi k\Bigr)
\\+\varphi(t+2\pi k)\frac d{dt}\mbox{sinc}\Bigl(\frac t2+\pi k\Bigr)\Biggr\}
e^{-i(t+2\pi k)\delta},
\end{multline}
while the second-order derivative is given by
\begin{widetext}
\begin{multline}
\frac{d^2\Phi_t}{dt^2}
=\sum_{k=-\infty}^\infty\Biggl\{
\left[\varphi''(t+2\pi k)-2i\delta\varphi'(t+2\pi k)
-\delta^2\varphi(t+2\pi k)\right]
\,\mbox{sinc}\Bigl(\frac t2+\pi k\Bigr)
\\
+2\left[\varphi'(t+2\pi k)-i\delta\varphi(t+2\pi k)\right]
\frac d{dt}\mbox{sinc}\Bigl(\frac t2+\pi k\Bigr)
+\varphi(t+2\pi k)\frac{d^2}{dt^2}\mbox{sinc}\Bigl(\frac t2+\pi k\Bigr)\Biggr\}
e^{-i(t+2\pi k)\delta}.
\end{multline}
\end{widetext}
Here for the sake of clarity, we have not yet introduced the explicit
forms for the derivatives of the sinc function,
\begin{equation}\begin{split}
\frac d{dt}\,\mbox{sinc}\Bigl(\frac{t}2+\pi k\Bigr)
&=\frac12\frac{\cos(t/2+\pi k)-\mbox{sinc}(t/2+\pi k)}{t/2+\pi k}
\\
\frac{d^2}{dt^2}\,\mbox{sinc}\Bigl(\frac{t}2+\pi k\Bigr)
&=-\frac14\mbox{sinc}\Bigl(\frac{t}2+\pi k\Bigr)
\\&-\frac12\frac{\cos(t/2+\pi k)-\mbox{sinc}(t/2+\pi k)}{(t/2+\pi k)^2}.
\end{split}\end{equation}
Next given that
\begin{equation}
\frac{d^k\Phi_t}{dt^k}=i^k\sum_{n=-\infty}^\infty\!n^ke^{itn}F_n\
\Rightarrow\
\sum_{n=-\infty}^\infty\!n^kF_n
=\frac1{i^k}\left.\frac{d^k\Phi_t}{dt^k}\right\rvert_{t=0},
\end{equation}
we can find that
\begin{equation}\label{eq:res1}\begin{split}
m&=\sum_{n=-\infty}^\infty\!nF_n=\frac1i\left.\frac{d\Phi_t}{dt}\right\rvert_{t=0}
=\mu-\delta+S_0
\\s^2&=\sum_{n=-\infty}^\infty\!(n-m)^2F_n
=-\left.\frac{d^2\Phi_t}{dt^2}\right\rvert_{t=0}-m^2
\\&=\sigma^2+\frac1{12}-S_1-S_0^2,
\end{split}\end{equation}
where
\begin{equation}\label{eq:msd}\begin{split}
S_0&=\sum_{\substack{k=-\infty\\k\ne0}}^\infty\!
\frac{(-1)^k\varphi(2\pi k)}{2\pi ik}e^{-2\pi ik\delta}
\\&=\sum_{\substack{k=-\infty\\k\ne0}}^\infty\!
\frac{(-1)^k\tilde\varphi(2\pi k)}{2\pi ik}e^{2\pi ik(\mu-\delta)},
\end{split}\end{equation}
and
\begin{equation}\label{eq:ssd}\begin{split}
S_1&=\sum_{\substack{k=-\infty\\k\ne0}}^\infty\!\frac{(-1)^k}{\pi k}
\left[\varphi'(2\pi k)-i\mu\varphi(2\pi k)-\frac{\varphi(2\pi k)}{2\pi k}\right]
e^{-2\pi ik\delta}
\\&=\sum_{\substack{k=-\infty\\k\ne0}}^\infty\!\frac{(-1)^k}{\pi k}
\left[\tilde\varphi'(2\pi k)-\frac{\tilde\varphi(2\pi k)}{2\pi k}\right]
e^{2\pi ik(\mu-\delta)}.
\end{split}\end{equation}
Here we have used the fact that $\mbox{sinc}(\pi k)=(\pi k)^{-1}\sin(\pi k)=0$
for any non-zero integer $k\in\mathbb Z-\{0\}$ as well as
\begin{equation}
\left.\frac d{dt}\,\mbox{sinc}\Bigl(\frac{t}2+\pi k\Bigr)\right\rvert_{t=0}
=\begin{cases}
\dfrac{(-1)^k}{2\pi k}&k\in\mathbb Z-\{0\}
\smallskip\\0&k=0\end{cases},
\end{equation}
and
\begin{equation}
\left.\frac{d^2}{dt^2}\,\mbox{sinc}\Bigl(\frac{t}2+\pi k\Bigr)\right\rvert_{t=0}
=\begin{cases}
\dfrac{(-1)^{k+1}}{2(\pi k)^2}&k\in\mathbb Z-\{0\}
\smallskip\\-\dfrac1{12}&k=0\end{cases}.
\end{equation}
Equations (\ref{eq:res1}) indeed reproduce
the results expected from the elementary arguments given
in the introduction with the proviso that the infinite sums,
$S_0$ and $S_1$ in equations (\ref{eq:msd}) and (\ref{eq:ssd})
are negligible. In other words, if one considers only the $k=0$
term in the characteristic function of equation (\ref{eq:fn_char}),
we would recover the results that $m=\mu-\delta$ and $s^2=\sigma^2+1/12$.
In fact, if one regards $F_n$ to be a continuous distribution over
real $n$ and replace the infinite sum $\sum_{n=-\infty}^\infty e^{i(t-\tau)n}$
in equation (\ref{eq:fn_char}) with the integral
$\int_{-\infty}^\infty e^{i(t-\tau)n}{\rm d}n$,
the Dirac comb $\sum_{k=-\infty}^\infty\deltaup(t-\tau+2\pi k)$
would be replaced by a single Dirac delta $\deltaup(t-\tau)$.
That is to say, the naive expectation that
$m=\mu-\delta$ and $s^2=\sigma^2+1/12$ may be considered as
the approximation in the limit of continuous $F_n$.

\subsection{symmetric distribution}

If $M\in\mathbb Z$ is the integer to which the mean $\mu$ is rounded,
$\mu\in[M+\delta-1/2,M+\delta+1/2]$ and $\chi=\mu-\delta-M\in[-1/2,1/2]$.
Since $M$ and $k$ are integers and $\mu-\delta=M+\chi$, we find
\begin{equation}\label{eq:csum}\begin{split}
S_0&=\sum_{\substack{k=-\infty\\k\ne0}}^\infty\!
\frac{(-1)^k\tilde\varphi(2\pi k)}{2\pi ik}e^{2\pi ik\chi};
\\S_1&=\sum_{\substack{k=-\infty\\k\ne0}}^\infty\!\frac{(-1)^k}{\pi k}
\left[\tilde\varphi'(2\pi k)-\frac{\tilde\varphi(2\pi k)}{2\pi k}\right]
e^{2\pi ik\chi},
\end{split}\end{equation}
which is basically the Fourier series expressions
of $S_0(\chi)$ and $S_1(\chi)$ for $\chi\in[-1/2,1/2]$.
Since $\tilde\varphi(-t)=\tilde\varphi(t)$
and $\tilde\varphi'(-t)=-\tilde\varphi'(t)$
for a symmetric distribution,
these can be further reducible to the real ones:
\begin{equation}\label{eq:rsum}\begin{split}
S_0&=\sum_{k=1}^\infty\!
\frac{(-1)^k\tilde\varphi(2\pi k)}{\pi k}\sin(2\pi k\chi);
\\S_1&=\sum_{k=1}^\infty\!\frac{2\cdot(-1)^k}{\pi k}
\left[\tilde\varphi'(2\pi k)-\frac{\tilde\varphi(2\pi k)}{2\pi k}\right]
\cos(2\pi k\chi)
\\&=\sum_{k=1}^\infty\!4\cdot(-1)^k
\left.\frac d{dt}\biggl(\frac{\tilde\varphi(t)}t\biggr)\right|_{t=2\pi k}
\cos(2\pi k\chi)
\end{split}\end{equation}
if $f(x)$ is symmetric about its mean.

Next, consider the family of the distributions sharing the common normalized
form; namely,
\begin{equation}
f(x)=\frac1\sigma\,F\!\left(\frac{x-\mu}\sigma\right)
\end{equation}
where $F(u)$ is a fixed non-negative function such that
$\int_{-\infty}^\infty\!{\rm d}u\,F(u)=1$,
$\int_{-\infty}^\infty\!{\rm d}u\,u\,F(u)=0$, and
$\int_{-\infty}^\infty\!{\rm d}u\,u^2F(u)=1$.
Then, for all members of the family,
we find $\varphi(t)=e^{i\mu t}\Phi(\sigma t)$ and
$\tilde\varphi(t)=\Phi(\sigma t)$ where
\begin{equation}
\Phi(\tau)=\int_{-\infty}^\infty\!{\rm d}u\,e^{iu\tau}F(u)
\end{equation}
is the characteristic function of the normalized distribution.
Here $f(x)$ is a symmetric distribution if and only if
$F(u)$ and $\Phi(\tau)$ are even functions:
$F(-u)=F(u)$ and $\Phi(-\tau)=\Phi(\tau)$. In the limit of
$\sigma=0$ -- essentially $F(u)=\deltaup(u)$ -- we then have
$\tilde\varphi(t)=\Phi(0)=1$ and
$\tilde\varphi'(t)=\sigma\Phi'(\sigma t)=0$ and so
\begin{equation}\begin{split}
\lim_{\sigma^2\to0}S_0&=\sum_{k=1}^\infty\!
\frac{(-1)^k}{\pi k}\sin(2\pi k\chi)=-\chi;
\\\lim_{\sigma^2\to0}S_1&=\sum_{k=1}^\infty\!\frac{(-1)^{k+1}}{(\pi k)^2}
\cos(2\pi k\chi)=\frac1{12}-\chi^2,
\end{split}\end{equation}
for $\chi\in[-1/2,1/2]$. Then it follows that
$m=\mu-\delta-\chi=M$ and $s^2=\sigma^2+1/12-(1/12-\chi^2)-(-\chi)^2
=\sigma^2=0$,
as expected (i.e.\
every sample point is rounded to the same integer).

Now suppose that $\Phi(\tau)$ admits an asymptotic expansion;
\begin{equation}\label{eq:asymphi}
\Phi(\tau)\simeq\frac1{|\tau|^s}
\sum_{p=0}^\infty\frac{\Phi_{\infty,p}}{\tau^{2p}}
\qquad\text{($\tau\to\infty$)}
\end{equation}
with a constant $s>0$. Then it follows from equation (\ref{eq:rsum}) that
\begin{equation}\label{eq:sasum}\begin{split}
S_0&\simeq\sum_{p=0}^\infty\!\frac{\Phi_{\infty,p}}{(2\sigma)^{2p+s}}
\sum_{k=1}^\infty\!\frac{(-1)^k\sin(2\pi k\chi)}{(\pi k)^{2p+s+1}};
\\S_1&\simeq\sum_{p=0}^\infty\!
\frac{(2p+s+1)\Phi_{\infty,p}}{(2\sigma)^{2p+s}}
\sum_{k=1}^\infty\!\frac{(-1)^{k+1}\cos(2\pi k\chi)}{(\pi k)^{2p+s+2}}.
\end{split}\end{equation}
Here the inner sums on $k$ converge absolutely for $s>0$ (NB: the sums
for an even integer $s$ are actually reducible to the Bernoulli polynomials)
given that
\begin{multline}
\left\lvert
\sum_{k=1}^\infty\!\frac{(-1)^k\sin(2\pi k\chi)}{(\pi k)^{2p+s+1}}
\right\rvert\le
\sum_{k=1}^\infty\!\frac{|(-1)^k\sin(2\pi k\chi)|}{(\pi k)^{2p+s+1}}
\\\le\sum_{k=1}^\infty\!\frac1{(\pi k)^{2p+s+1}}
=\frac{\zeta(2p+s+1)}{\pi^{2p+s+1}},
\end{multline}
where $\zeta(z)$ is the Riemann zeta function and similarly
\begin{equation}
\left\lvert
\sum_{k=1}^\infty\!\frac{(-1)^{k+1}\cos(2\pi k\chi)}{(\pi k)^{2p+s+2}}
\right\rvert\le\frac{\zeta(2p+s+2)}{\pi^{2p+s+2}}.
\end{equation}
If $\chi=1/2$, then $\cos(2\pi k\chi)=\cos(\pi k)=(-1)^k$
for any integer $k$ and so the last bound is actually sharp.
By contrast, the first bound is not sharp but it suffices for our purposes here.
Since $\lim_{z\to\infty}\zeta(z)=1$ and $\zeta(z)$ for $z>1$
is monotonically decreasing, we can conclude that
equations (\ref{eq:sasum}) is in fact valid asymptotic
expansion of $S_0$ and $S_1$ as $\sigma\to\infty$.
Also it follows that, if
$\lim_{\tau\to\infty}d\ln|\Phi(\tau)|/d\ln|\tau|=-s<0$,
we have $S_0\sim\sigma^{-s}\to0$ and $S_1\sim\sigma^{-s}\to0$
as $\sigma\to\infty$ as well as $m=\mu-\delta+\mathcal O(\sigma^{-s})$
and $s^2=\sigma^2+1/12+\mathcal O(\sigma^{-s})$.

As a concrete example, consider the bilateral exponential
(Laplace) distribution of the mean $\mu$ and the variance $\sigma^2$:
\begin{equation}\label{eq:bed}
f(x)=\frac1{\!\sqrt2\sigma}
\exp\biggl(-\frac{\!\sqrt2}\sigma\lvert x-\mu\rvert\biggr),
\end{equation}
which is easily normalizable so that
\begin{equation}
F(u)=\frac{e^{-\!\sqrt2|u|}}{\!\sqrt2};\quad
\Phi(\tau)=\left(1+\frac{\tau^2}2\right)^{-1}.
\end{equation}
Here we find $\Phi(\tau)=-\sum_{k=1}^\infty(-2/\tau^2)^k\simeq2/\tau^2$
as $\tau\to\infty$ and so it should be $m=\mu-\delta+\mathcal O(\sigma^{-2})$
and $s^2=\sigma^2+1/12+\mathcal O(\sigma^{-2})$. In fact for this case,
we know the analytic forms for $S_0(\chi)$ and $S_1(\chi)$.
That is to say, let us consider the odd function for $\chi\in[-1/2,1/2]$
given by
\begin{equation}\label{eq:s0_bed}\begin{split}
S_0(\chi)&=\sum_{k=1}^\infty\frac{2^kB_{2k+1}(1/2+\chi)}{(2k+1)!\sigma^{2k}}
\\&=\frac{\sinh(\!\sqrt2\chi/\sigma)}{2\sinh[1/(\!\sqrt2\sigma)]}-\chi,
\end{split}\end{equation}
where $B_n(z)$ is the Bernoulli Polynomial. Then
we find that
\begin{equation}
\int_{-1/2}^{1/2}\!{\rm d}\chi\,S_0(\chi)\sin(2\pi k\chi)
=\frac{(-1)^k}{2\pi k[1+2(\sigma\pi k)^2]}.
\end{equation}
for a positive integer $k$.
It follows that the first infinite sum in equation (\ref{eq:rsum})
with $\tilde\varphi(t)=\Phi(\sigma t)=[1+(\sigma t)^2/2]^{-1}$ is
the Fourier (sine) series expansion for $S_0(\chi)$
in equation (\ref{eq:s0_bed}). In other words, if we sample
the random variable $x$ distributed according to equation (\ref{eq:bed})
and round it to an integer $n$ such that
$n+\delta-1/2\le x<n+\delta+1/2$ or $n+\delta-1/2<x\le n+\delta+1/2$
with a fixed $\delta\in[-1/2,1/2]$, the mean of $n$ is
\begin{equation}\begin{split}
m&=\mu-\delta+S_0
=M+\frac{\sinh(\!\sqrt2\chi/\sigma)}{2\sinh[1/(\!\sqrt2\sigma)]}
\\&\simeq\mu-\delta-\frac{\chi(1-4\chi^2)}{12\sigma^2}+\mathcal O(\sigma^{-4}),
\end{split}\end{equation}
where $M$ is the integer to which $\mu$ is rounded and
$\chi=\mu-\delta-M$. Similarly we can also establish that
the second infinite sum in equation (\ref{eq:rsum})
with the same $\tilde\varphi(t)$ is a Fourier series representation of
\begin{equation}\begin{split}
S_1(\chi)&=\sum_{k=1}^\infty\!\frac{2^kB_{2k+2}(1/2+\chi)}{(k+1)(2k)!\sigma^{2k}}
\\&=\sigma^2+\frac1{12}-\chi^2
+\frac{\chi\sinh(\!\sqrt2\chi/\sigma)}{\sinh[1/(\!\sqrt2\sigma)]}
\\&\phantom{=\sigma^2}
-\frac{\cosh[1/(\!\sqrt2\sigma)]\cosh(\!\sqrt2\chi/\sigma)}
{2\sinh^2[1/(\!\sqrt2\sigma)]},
\end{split}\end{equation}
and so the variance of the rounded integers sampled
over the random variables with the distribution
in equation (\ref{eq:bed}) is
\begin{equation}\label{eq:bedsasym}\begin{split}
s^2&=\frac{2\cosh[1/(\!\sqrt2\sigma)]\cosh(\!\sqrt2\chi/\sigma)-\sinh^2(\!\sqrt2\chi/\sigma)}{4\sinh^2[1/(\!\sqrt2\sigma)]}
\\&\simeq\sigma^2+\frac1{12}
-\left(\frac7{480}-\frac{\chi^2}4
+\frac{\chi^4}2\right)\frac1{\sigma^2}+\mathcal O(\sigma^{-4}).
\end{split}\end{equation}

\section{Expectation Values for unspecified mean}

The results so far have concerned the distributions with a known
fixed mean. Here instead we consider the cases of unspecified means.
That is to say, let us calculate the expectation values for
the (difference to the true) mean and the variance of the rounded
variables averaged over distributions with all possible means.
In practice, this is achieved by averaging over $\chi\in[-1/2,1/2]$
and so $\langle m-\mu\rangle=-\delta+\langle S_0\rangle$
and $\langle s^2\rangle=\sigma^2+1/12
-\langle S_1\rangle-\langle S_0^2\rangle$, where
$\langle S_0\rangle=\int_{-1/2}^{1/2}{\rm d}\chi\,S_0(\chi)$
and so on. However, we have
$\langle e^{2\pi ik\chi}\rangle
=\langle\sin(2\pi\chi k)\rangle=\langle\cos(2\pi\chi k)\rangle=0$
for any non-zero integer $k$ and thus $\langle S_0\rangle=\langle S_1\rangle=0$
given equations (\ref{eq:csum}) and (\ref{eq:rsum}).
As for $\langle S_0^2\rangle$, let us first note
$\langle S_0^2\rangle\ne\langle S_0\rangle^2$. Rather
from equation (\ref{eq:csum}),
\begin{equation}\label{eq:avs0}\begin{split}
\langle S_0^2\rangle&=
\sum_{\substack{k,p=-\infty\\k,p\ne0}}^\infty\!(-1)^{k+p}
\frac{\tilde\varphi(2\pi k)\tilde\varphi(2\pi p)}{(2\pi i)^2kp}
\langle e^{2\pi i(k+p)\chi}\rangle
\\&=\sum_{\substack{k=-\infty\\k\ne0}}^\infty\!
\frac{\tilde\varphi(2\pi k)\tilde\varphi(-2\pi k)}{(2\pi k)^2}
=\sum_{k=1}^\infty\!
\frac{\lvert\tilde\varphi(2\pi k)\rvert^2}{2(\pi k)^2},
\end{split}\end{equation}
where we have used the fact that $\tilde\varphi(-t)$ is
the complex conjugate of $\tilde\varphi(t)$ for any real $f(x)$
(see eq.~\ref{eq:scfdef}). The same result for the symmetric distributions
may also be derived from equation (\ref{eq:rsum}) given
$\langle\sin^2(2\pi\chi k)\rangle=1/2$ 
and $\langle\sin(2\pi\chi k)\sin(2\pi\chi p)\rangle=0$ for
positive integers $k\ne p$. Consequently
\begin{equation}\label{eq:s2av}
\langle m\rangle=\mu-\delta,\quad
\langle s^2\rangle=\sigma^2+\frac1{12}-\sum_{k=1}^\infty\!
\frac{\lvert\tilde\varphi(2\pi k)\rvert^2}{2(\pi k)^2},
\end{equation}
with $\lim_{\sigma\to0}\langle S_0^2\rangle=\zeta(2)/(2\pi^2)=1/12$
(given $\tilde\varphi(t)=1$ for $\sigma=0$) and
$\langle s^2\rangle=\sigma^2=0$ in the limit of $\sigma=0$.

If $\Phi(\tau)$ is given by the same function admitting
the asymptotic expansion of equation (\ref{eq:asymphi}),
we find
\begin{equation}
\langle S_0^2\rangle\simeq\sum_{p=0}^\infty\frac{\zeta(2p+2s+2)
\left(\sum_{r=0}^p\Phi_{\infty,p-r}\Phi_{\infty,r}\right)}
{2^{2(p+s)+1}\pi^{2(p+s+1)}\sigma^{2(p+s)}},
\end{equation}
and so $\langle S_0^2\rangle\sim\sigma^{-2s}\to0$
and $\langle s^2\rangle=\sigma^2+1/12+\mathcal O(\sigma^{-2s})$
(also $0\le\langle s^2\rangle\le\sigma^2+1/12$)
as $\sigma\to\infty$ for $\Phi(\tau)\sim\tau^{-s}$ as $\tau\to\infty$.
That is to say, $\langle s^2\rangle$ typically tends to the limiting value
$\lim_{\sigma^2\to\infty}(\langle s^2\rangle-\sigma^2)=1/12$
about twice much faster than the individual $s^2$ does.
For example, with the bilateral exponential distribution given
in equation (\ref{eq:bed}), we specifically have
\begin{equation}\begin{split}
\langle S_0^2\rangle&=\sum_{k=1}^\infty\!\frac1{2(\pi k)^2[1+2(\pi k\sigma)^2]^2}
\\&=\sum_{k=1}^\infty\!\frac1{2^3(\pi k)^6\sigma^4}
\sum_{p=0}^\infty\!\frac{(-1)^p(p+1)}{2^p(\pi k\sigma)^{2p}}
\\&=\sum_{p=0}^\infty\!\frac{(-1)^p(p+1)}{\sigma^{2(p+2)}}
\frac{\zeta(2p+6)}{(2\pi^2)^{p+3}}
\\&=\sum_{p=0}^\infty\!\frac{(p+1)B_{2p+6}}{(2p+6)!}\frac{2^{p+2}}{\sigma^{2(p+2)}}
\\&=\sigma^2+\frac1{12}
-\frac{2+3\!\sqrt2\,\sigma\sinh(\!\sqrt2/\sigma)}{16\sinh^2[1/(\!\sqrt2\sigma)]}
\end{split}\end{equation}
where $B_k$ is the Bernoulli number, and so
\begin{equation}\begin{split}
\langle s^2\rangle
&=\frac{2+3\!\sqrt2\,\sigma\sinh(\!\sqrt2/\sigma)}{16\sinh^2[1/(\!\sqrt2\sigma)]}
\\&=\sigma^2+\frac1{12}
-\frac1{\sigma^4}
\sum_{p=0}^\infty\!\frac{(p+1)B_{2p+6}}{(2p+6)!}\frac{2^{p+2}}{\sigma^{2p}};
\end{split}\end{equation}
that is,
$\langle s^2\rangle\simeq\sigma^2+1/12-1/(7560\sigma^4)
+\mathcal O(\sigma^{-6})$, which contrasts to equation (\ref{eq:bedsasym}).

\subsection{distributions with a compact support}

Suppose that $f(x)$ is 
\begin{equation}\label{eq:ud}
f(x)=\begin{cases}
\dfrac1{2\!\sqrt3\sigma}&(\mu-\!\sqrt3\sigma\le x\le\mu+\!\sqrt3\sigma)
\\0&\text{elsewhere}\end{cases},
\end{equation}
i.e.\ the uniform distribution over a compact interval,
the normalized form of which is
\begin{equation}\begin{split}
F(u)&=\begin{cases}
\dfrac1{2\!\sqrt3}&(-\!\sqrt3\le u\le\!\sqrt3)
\\0&\text{elsewhere}\end{cases};
\\\Phi(\tau)&=\int_{-\!\sqrt3}^{\!\sqrt3}
\frac{e^{iu\tau}{\rm d}u}{2\!\sqrt3}
=\frac{\sin(\!\sqrt 3\tau)}{\!\sqrt3\tau}.
\end{split}\end{equation}
We then find for the compact uniform distribution that
\begin{equation}\label{eq:s0unid}\begin{split}
\langle S_0^2\rangle&=\sum_{k=1}^\infty\!
\frac{\sin^2(2\!\sqrt3\pi\sigma k)}{2(\pi k)^2(2\!\sqrt3\pi\sigma k)^2}
=\sum_{k=1}^\infty\!\frac{1-\cos(4\!\sqrt3\pi\sigma k)}{48\sigma^2(\pi k)^4}
\\&=\frac1{48\sigma^2}\left[
\frac{\zeta(4)}{\pi^4}-\frac{B_4(\xi)}3\right]
=\frac{\xi^2(1-\xi)^2}{144\sigma^2}
\end{split}\end{equation}
where $\xi=2\!\sqrt3\sigma-\lfloor2\!\sqrt3\sigma\rfloor\in[0,1)$
is the fractional part of ``$2\!\sqrt3\sigma$
(which is the width of the support)''. Here we have used
the Fourier series expansion of the Bernoulli polynomial
(for $0\le\xi\le1$) of the even order\citep{DLMF}
\begin{equation}
B_{2n}(\xi)=(-1)^{n+1}\frac{(2n)!}{2^{2n-1}}
\sum_{k=1}^\infty\frac{\cos(2\pi k\xi)}{(\pi k)^{2n}},
\end{equation}
with $\zeta(4)=\pi^4/90$ and $B_4(\xi)=\xi^2(\xi-1)^2-1/30$.
That is to say, while the remainder
$\langle S_0^2\rangle=\sigma^2+1/12-\langle s^2\rangle$ falls off
``on average'' like $\sim\sigma^{-2}$ as $\sigma\to\infty$,
the actual behavior includes the periodic modulation
superimposed on the asymptotic scale-free decay. This is due
to the compact support on the underlying distribution: note
that a compact distribution $F(u)$ typically results in
an oscillatory $\Phi(\tau)$, and $\langle S_0^2\rangle$ is basically
the sum on the regular sampling of the latter. The resulting
modulation of $\langle S_0^2\rangle$ may be understood as a sort
of interference patterns between the width of the compact support
and the unit intervals for the rounded integer values.
However unless the variance of the underlying continuous
distribution $f(x)$ is known a priori, the averaged asymptotic
behavior of $\langle S_0^2\rangle$ as $\sigma\to\infty$ may
be used to estimate $\sigma^2$ from $s^2$ in practice within
a reasonable accuracy (provided $s^2>>1$). If one is in fact only
interested in the averaged asymptotic behavior, we can further
average $\langle S_0^2\rangle$ in equation (\ref{eq:s0unid}) over
$\xi\in[0,1)$ and get
$\langle S_0^2\rangle=\zeta(4)/(48\sigma^2\pi^4)=(4320\sigma^2)^{-1}$.
For a general distribution with a compact support, one may obtain
the averaged asymptotic behavior for $\langle S_0^2\rangle$ in
equation (\ref{eq:avs0}) by assuming that any sum of the form
$\sum_k\sin(a\sigma k)/k^n$ or $\sum_k\cos(a\sigma k)/k^n$
(where $a$ is a fixed real constant) also vanishes on average.

In principle we can also calculate $S_0$ and $S_1$ first,
and subsequently average them over the proper interval.
For the uniform distribution in equation (\ref{eq:ud}),
equation (\ref{eq:rsum}) results in
\begin{equation}\label{eq:us0}\begin{split}
S_0&=\sum_{k=1}^\infty\!(-1)^k
\frac{\sin(2\!\sqrt3\sigma\pi k)}{2\!\sqrt3\sigma(\pi k)^2}\sin(2\pi k\chi)
\\&=\sum_{k=1}^\infty\!
\frac{\cos(2\pi\Delta_-k)-\cos(2\pi\Delta_+k)}{4\!\sqrt3\sigma(\pi k)^2}
\\&=\frac{B_2(\Delta_-)-B_2(\Delta_+)}{4\!\sqrt3\sigma}
=-\frac{\lambda\zeta}{\!\sqrt3\sigma},
\end{split}\end{equation}
where $\Delta_\pm=\mu-\delta+1/2\pm\!\sqrt3\sigma-m_\pm$ is
the fractional part of $\mu\pm\!\sqrt3\sigma-\delta+1/2$
and $m_\pm=\lfloor(1/2+\mu-\delta\pm\!\sqrt3\sigma)\rfloor$
is the integer to which the upper/lower limit
of the compact support (i.e.\ $\mu\pm\!\sqrt3\sigma$) is rounded.
In addition, $\lambda=(\Delta_++\Delta_--1)/2$ and
$\zeta=(\Delta_+-\Delta_-)/2$. Also used are
$(-1)^k\sin(2\pi k\chi)=\sin[2\pi k(1/2+\mu-\delta)]$
for any integer $k$ given $\chi=\mu-\delta-M$ with an integer $M$,
and $B_2(x)=x^2-x+1/6$. Similarly,
\begin{equation}\label{eq:us1}\begin{split}
S_1&=\sum_{k=1}^\infty\!(-1)^k
\left[\frac{\cos(2\!\sqrt3\sigma\pi k)}{(\pi k)^2}
-\frac{\sin(2\!\sqrt3\sigma\pi k)}{\!\sqrt3\sigma(\pi k)^3}\right]
\cos(2\pi k\chi)
\\&=\frac{B_2(\Delta_-)+B_2(\Delta_+)}{2}
+\frac{B_3(\Delta_-)-B_3(\Delta_+)}{3\!\sqrt3\sigma}
\\&=\lambda^2+\zeta^2-\frac1{12}
-\left(2\lambda^2+\frac23\zeta^2-\frac16\right)\frac\zeta{\!\sqrt3\sigma},
\end{split}\end{equation}
further utilizing $B_3(x)=x^3-3x^2/2+x/2$ and
the Fourier series for the odd-order Bernoulli polynomial\citep{DLMF}:
\begin{equation}
\sum_{k=1}^\infty\frac{\sin(2\pi kx)}{(\pi k)^{2n+1}}
=\frac{(-1)^{n+1}2^{2n}}{(2n+1)!}B_{2n+1}(x)
\quad(0\le x\le 1).
\end{equation}
Here note $\lim_{\sigma^2\to0}S_1=\lambda^2+\zeta^2-1/12$ and so
$s^2=\sigma^2-\lambda^2-\zeta^2+1/6+\mathcal O(\sigma^{-1})$
even though $\Phi(\tau)=\mbox{sinc}(\!\sqrt3\tau)\sim\tau^{-1}$
as $\tau\to\infty$. Technically this is not a counter-example of
the prior discussion following equation (\ref{eq:asymphi}), since
$\mbox{sinc}(x)$ does not actually have a proper asymptotic
expansion as $x\to\infty$ in the strict sense. In fact,
we observe that, while the asymptotic behavior of $S_0$
as $\sigma\to\infty$ follows that of $\Phi(\tau)$ as
$\tau\to\infty$, the behavior of $S_1$ actually traces
$\tau\Phi'(\tau)$ instead. With an oscillatory $\Phi(\tau)$
due to $F(u)$ on a compact support, $\tau\Phi'(\tau)$ can indeed
be much larger than $\Phi(\tau)$ even if $\Phi(\tau)$ is bounded by
an asymptotically decaying envelope, and so $s^2-\sigma^2$ does
not necessarily tend to $1/12$ unless $\tau\Phi'(\tau)\to0$ as $\tau\to\infty$.
By contrast, the formula for $\langle s^2\rangle$ in equation (\ref{eq:s2av})
only involves $|\tilde\varphi(2\pi k)|^2$ and
we thus expect the asymptotic behavior of the remainder
$\sigma^2+1/12-\langle s^2\rangle$ to trace that of $|\Phi(\tau)|^2$ in general.

In order to average the expresseions
in equations (\ref{eq:us0}) and (\ref{eq:us1}) over $\chi\in[-1/2,1/2]$,
we first observe that $\zeta=\xi/2\ge0$ or $\zeta=(\xi-1)/2<0$
(where $0\le\xi<1$ is still the fractional part of $2\!\sqrt3\sigma$).
Next $\Delta_\pm=\lambda\pm\zeta+1/2\in[0,1)$ implies
$\lambda\in[-1/2+\zeta,1/2-\zeta)=[(\xi-1)/2,(1-\xi)/2)$ for
$\zeta\ge0$ and
$\lambda\in[-1/2-\zeta,1/2+\zeta)=[-\xi/2,\xi/2)$ for $\zeta<0$.
Finally $\lambda=\chi\pm\frac12$ or $\chi$ depending on the parity
of $m_++m_-$ and we find that $\chi\in[-1/2,1/2]$ at fixed
$\sigma$ (and consequently $\xi$ is also fixed) then maps to the union of
those two allowed intervals on $\lambda$. Consequently
averaging over $\chi\in[-1/2,1/2]$ is achieved through 
summing two integrals on $\lambda$: viz.\
\begin{equation}\label{eq:s0av}\begin{split}
\langle S_0^2\rangle
&=\int_{\frac{\xi-1}2}^{\frac{1-\xi}2}\!{\rm d}\lambda
\left.S_0^2\right\rvert_{\zeta=\frac\xi2}
+\int_{-\frac\xi2}^{\frac\xi2}\!{\rm d}\lambda
\left.S_0^2\right\rvert_{\zeta=\frac{\xi-1}2}
\\&=\frac2{3\sigma^2}
\left[\frac{\xi^2}4\int_0^{\frac{1-\xi}2}\!{\rm d}\lambda\,\lambda^2
+\frac{(\xi-1)^2}4\int_0^{\frac\xi2}\!{\rm d}\lambda\,\lambda^2\right]
\\&=\frac2{3\sigma^2}
\left[\frac{\xi^2}4\frac{(1-\xi)^3}{24}
+\frac{(\xi-1)^2}4\frac{\xi^3}{24}\right]
=\frac{\xi^2(1-\xi)^2}{144\sigma^2},
\end{split}\end{equation}
which recovers equation (\ref{eq:s0unid}). Since $S_0$ is an odd function
of $\lambda$, it is immediately obvious that $\langle S_0\rangle=0$.
As for $\langle S_1\rangle$, two respective integrals for $\zeta=\xi/2\ge0$
and $\zeta=(\zeta-1)/2<0$ exactly cancel each other and so
$\langle S_1\rangle=0$.

Since $\chi$ (i.e.\ the offset of the mean from an integer value)
determines both $\Delta_\pm$ (i.e.\ the offsets of the boundary
points of the support from integer values) once $\sigma$ (which
specifies the width of the support) is fixed, $\Delta_+$ and $\Delta_-$
are not in fact independent from each other. Nevertheless it is still
formally possible to consider $S_0$ in equation (\ref{eq:us0}) as
a function of the pair of independent variables
$(\Delta_+,\Delta_-)\in[0,1)^2$ and average $S_0^2$
over this whole rectangular domain. Following the coordinate transform
$(\Delta_+,\Delta_-)\rightarrow(\lambda,\zeta)$,
the resulting average is shown to be identical to further
averaging the $\chi$-average of equation (\ref{eq:s0av})
over $\zeta\in(-1/2,1/2)$ or equivalently $\xi\in[0,1)$.
That is to say, if one assume that both upper and lower boundaries of
the compact support are randomly placed relative to the integer values
(and independent from each other), the resulting expectation value
(averaged over all possible such placements) recovers only
the ``slow'' asymptotic decay behavior of $\langle s^2\rangle$
for the rounded random variables on a compact support while averaging
off the ``fast'' modulation due to the interference between the width
of the support and the integer signposts. 

\section{Normal distributions}

Finally we would like to consider the case of $f(x)$ being
the Gaussian normal distribution:
\begin{equation}
f(x)=\frac1{(2\pi)^{1/2}\sigma}
\exp\biggl[-\frac{(x-\mu)^2}{2\sigma^2}\biggr],
\end{equation}
or in the standard form
\begin{equation}
F(u)=\frac{e^{-u^2/2}}{\!\sqrt{2\pi}},\quad
\Phi(\tau)=e^{-\tau^2/2}.
\end{equation}
Then with $\tilde\varphi(t)=\Phi(\sigma t)$, we have
\begin{equation}\begin{split}
S_0&=\sum_{k=1}^\infty\!(-1)^k\sin(2\pi k\chi)\,
\frac{e^{-2(\pi\sigma k)^2}}{\pi k};
\\S_1&=\sum_{k=1}^\infty\!(-1)^k\frac{\cos(2\pi k\chi)}{\pi^2}
\,\left.
\frac d{d\kappa}\biggl(\frac{e^{-2(\pi\sigma\kappa)^2}}\kappa\biggr)
\right\rvert_{\kappa=k}
\\&=\sum_{k=1}^\infty\!(-1)^{k+1}\cos(2\pi k\chi)\,
\frac{1+(2\pi\sigma k)^2}{(\pi k)^2}e^{-2(\pi\sigma k)^2}.
\end{split}\end{equation}
Thanks to the super-exponential decay $\propto e^{-2(\pi\sigma k)^2}$ in $k$
(NB: the $k=2$ term is suppressed relative to the $k=1$ term
by $\sim e^{-6(\pi\sigma)^2}$: if $\sigma=1$,
note $e^{-6\pi^2}\approx2\times10^{-26}$!),
these sums (which are actually in the form of the \emph{Jacobi theta function}
and its antiderivatives) converge extremely quickly and are completely
dominated by their respective first terms unless $\sigma\ll1$.
Alternatively we may construct a more formal argument
by bracketing the infinite sums following the integral convergence test.
In particular, we first observe that
\begin{equation}\begin{split}
|S_0|&\le\sum_{k=1}^\infty\frac{e^{-2(\pi\sigma k)^2}}{\pi k}=I_0,
\\|S_1|&\le\sum_{k=1}^\infty
\frac{1+(2\pi\sigma k)^2}{(\pi k)^2}e^{-2(\pi\sigma k)^2}=I_1
\end{split}\end{equation}
but the summands now are strictly decreasing positive functions
of $k\ge1$. Hence the integral convergence test indicates
\begin{equation}\begin{split}
\frac{E_1(2\pi^2\sigma^2)}{2\pi}&\le
I_0\le\frac{e^{-2(\pi\sigma)^2}}\pi+\frac{E_1(2\pi^2\sigma^2)}{2\pi},
\\\frac{e^{-2(\pi\sigma)^2}}{\pi^2}&\le
I_1\le\frac{(2\pi\sigma)^2+2}{\pi^2}e^{-2(\pi\sigma)^2},
\end{split}\end{equation}
where $E_1(x)=\int_1^\infty\!{\rm d}t\,e^{-tx}/t$
is the analytic exponential integral. Note that $I_0$ for $\sigma=0$,
which results in the harmonic series, actually diverges
and so the first bounds are only valid for $\sigma>0$,
but it has been already shown that $S_0=-\chi$ if $\sigma=0$.
Given the asymptotic expansion $e^xE_1(x)\sim\sum_{k=0}^\infty(-1)^kk!/x^{k+1}$
as $x\to\infty$, the first bounds may also be replaced by the purely
elementary functions. In particular, for $x>0$, we find
\begin{equation}
e^xE_1(x)=\int_1^\infty\!\frac{e^{x(1-t)}{\rm d}t}t
<\int_1^\infty\!{\rm d}t\,e^{-x(t-1)}=\frac1x,
\end{equation}
i.e.\ $E_1(x)<e^{-x}/x$ for $x>0$, and so follows that
\begin{equation}
|S_0|\le I_0\le
\left[1+\frac1{(2\pi\sigma)^2}\right]\frac{e^{-2(\pi\sigma)^2}}\pi.
\end{equation}
That is to say, for a sufficiently large $\sigma$,
both sums $I_0$ and $I_1$ are completely dominated
by their respective first terms, and $S_0\sim\mathcal O(e^{-2(\pi\sigma)^2})$
and $S_1\sim\mathcal O(\sigma^2e^{-2(\pi\sigma)^2})$
as $\sigma\to\infty$. In conclusion, the mean $m$ and the variance $s^2$
of the rounded variables drawn from the normal distribution (of
the mean $\mu$ and the variance $\sigma^2$) behave like
\begin{equation}\begin{split}
m&\simeq\mu-\delta+\mathcal O(e^{-2(\pi\sigma)^2});
\\s^2&\simeq\sigma^2+\frac1{12}+\mathcal O(\sigma^2e^{-2(\pi\sigma)^2}),
\end{split}\end{equation}
as $\sigma\to\infty$, but the remainder in most practical purposes
can be safely ignored provided $\sigma\gtrsim1$
(NB: $e^{-2\pi^2}/\pi\approx8.5\times10^{-10}$).

As for the expectation value averaged over $\mu$ at fixed $\sigma^2$,
if we consider the sum
\begin{equation}
\langle S_0^2\rangle=
\sum_{k=1}^\infty\frac{e^{-(2\pi\sigma k)^2}}{2(\pi k)^2},
\end{equation}
the integral test is still applicable:
\begin{equation}
\int_1^\infty\!\frac{e^{-(2\pi\sigma x)^2}{\rm d}x}{2(\pi x)^2}\le
\langle S_0^2\rangle\le
\frac{e^{-(2\pi\sigma)^2}}{2\pi^2}
+\int_1^\infty\!\frac{e^{-(2\pi\sigma x)^2}{\rm d}x}{2(\pi x)^2}.
\end{equation}
While here the integral is technically reducible to the incomplete
gamma function (or an expression involving the error function),
it is sufficient  for our purpose to note
\begin{equation}
e^a\!\int_1^\infty\!\frac{e^{-ax^2}{\rm d}x}{x^2}
=\int_0^\infty\!\frac{e^{-at}{\rm d}t}{2(t+1)^{3/2}}
<\int_0^\infty\!\frac{e^{-at}{\rm d}t}2=\frac1{2a}
\end{equation}
for $a>0$. Hence it follows that, for $\sigma>0$,
\begin{gather}
\langle S_0^2\rangle=\sigma^2+\frac12-\langle s^2\rangle
\le\left[\frac12+\frac1{(4\pi\sigma)^2}\right]\frac{e^{-(2\pi\sigma)^2}}{\pi^2};
\\\therefore\
\langle s^2\rangle\simeq\sigma^2+\frac1{12}+\mathcal O(e^{-(2\pi\sigma)^2})
\end{gather}
with an even faster-decaying (cf.\
$e^{-(2\pi)^2}/\pi^2\approx7.2\times10^{-19}$) remainder term.
In summary, the rounding errors for the normally distributed random variables
can for most practical applications be considered as independent from
the intrinsic dispersion unless the intrinsic dispersion itself is quite smaller
than the rounding unit.

\end{document}